\documentclass[a4paper,12pt,leqno]{article}
\usepackage{latexsym}

\usepackage{amssymb} 
\usepackage{amsmath} 
\usepackage{theorem}

\def\Z{{\mathbb{Z}}}
\def\R{{\mathbb{R}}}
\def\K{{\mathbb{K}}}
\def\A{{\mathcal{A}}}
\def\B{{\mathcal{B}}}

\DeclareMathOperator{\codim}{codim}

\DeclareMathOperator{\Der}{Der}

\numberwithin{equation}{section}

\newcommand{\owari}{\hfill$\square$}

\theoremstyle{break}
\newtheorem{theorem}{Theorem}[section]
\newtheorem{prop}[theorem]{Proposition}
\newtheorem{cor}[theorem]{Corollary}
\newtheorem{lemma}[theorem]{Lemma}
\newtheorem{define}[theorem]{Definition}

\newcommand{\xgraphAvertex}[1][****]{
\xgraphAVertex #1
}
\newcommand{\xgraphAVertex}[4]{
\if#3o\put(0,0){\circle{4}}\fi%
\if#3*\put(0,0){\circle*{4}}\fi%
\if#3.\put(0,0){\circle*{2.4}}\fi%
\if#4o\put(30,0){\circle{4}}\fi%
\if#4*\put(30,0){\circle*{4}}\fi%
\if#4.\put(30,0){\circle*{2.4}}\fi%
\if#2o\put(0,30){\circle{4}}\fi%
\if#2*\put(0,30){\circle*{4}}\fi%
\if#2.\put(0,30){\circle*{2.4}}\fi%
\if#1o\put(30,30){\circle{4}}\fi%
\if#1*\put(30,30){\circle*{4}}\fi%
\if#1.\put(30,30){\circle*{2.4}}\fi%
}
\newcommand{\xgraphA}[6]{
\if#1+\put(30,30){\line(-1,0){30}}\fi
\if#1.\qbezier[7](30,30)(15,30)(0,30)\fi
\if#1-\put(30,31){\line(-1,0){30}}\put(30,29){\line(-1,0){30}}\fi
\if#2+\put(0,0){\line(1,1){30}}\fi 
\if#2.\qbezier[10](0,0)(15,15)(30,30)\fi 
\if#2-\put(-0.7,0.7){\line(1,1){30}}\put(0.7,-0.7){\line(1,1){30}}\fi 
\if#3+\put(30,30){\line(0,-1){30}}\fi
\if#3.\qbezier[7](30,0)(30,15)(30,30)\fi 
\if#3-\put(29,30){\line(0,-1){30}}\put(31,30){\line(0,-1){30}}\fi
\if#4+\put(0,0){\line(0,1){30}}\fi 
\if#4.\qbezier[7](0,0)(0,15)(0,30)\fi 
\if#4-\put(-1,0){\line(0,1){30}}\put(1,0){\line(0,1){30}}\fi 
\if#5+\put(0,30){\line(1,-1){30}}\fi 
\if#5.\qbezier[10](30,0)(15,15)(0,30)\fi 
\if#5-\put(-0.7,29.3){\line(1,-1){30}}\put(0.7,30.7){\line(1,-1){30}}\fi 
\if#6+\put(0,0){\line(1,0){30}}\fi 
\if#6.\qbezier[7](0,0)(15,0)(30,0)\fi 
\if#6-\put(0,1){\line(1,0){30}}\put(0,-1){\line(1,0){30}}\fi 
\xgraphAvertex}

\title{Totally free arrangements of hyperplanes}
\author{Takuro Abe\thanks{
Department of Mathematics, Hokkaido University, 
Kita-10, Nishi-8, Kita-Ku, 
Sapporo, 060-0810, Japan.
email:abetaku@math.sci.hokudai.ac.jp.},
Hiroaki Terao\thanks{
Department of Mathematics, Hokkaido University, 
Kita-10, Nishi-8, Kita-Ku, 
Sapporo, 060-0810, Japan.
email:terao@math.sci.hokudai.ac.jp.} and
Masahiko Yoshinaga\thanks{
Department of Mathematics, Kobe University,
1-1 Rokkodai, Nada-ku,
Kobe, 657-8501, Japan.
email:myoshina@math.kobe-u.ac.jp.
}
}
\date{\today}

\begin{document}

\maketitle

\begin{abstract}
A
central arrangement $\A$ 
of hyperplanes 
in an $\ell$-dimensional  
vector space $V$ 
is said to be {\it
totally free}
if a multiarrangement $(\A, m)$ is free
for any multiplicity $ m : \A\rightarrow \Z_{> 0}$. 
It has been known that $\A$ 
is totally free whenever $\ell \le 2$.  
In this article, we will prove
that there does not exist any totally free
arrangement other than the obvious ones, that is, a product of 
one-dimensional arrangements and two-dimensional 
ones.
\end{abstract}

\setcounter{section}{0}
\section{Introduction}
Let $V$ be an $\ell$-dimensional
vector space $(\ell \ge 1)$ over $\K$ 
with a coordinate system $\{x_1,\ldots,x_\ell\} \subset V^*$.
Define 
$S:=\mbox{Sym}(V^*) \simeq \K[x_1,\ldots,x_\ell]$.
Let $\Der_{\K}(S) $ be the set of all $\K$-linear derivations of
$S$ to itself.  Then
$
\Der_\K(S)
=
\oplus_{i=1}^\ell S\cdot \partial_{x_i}$ is a free $S$-module of rank $\ell$. 
A \textit{central arrangement (of hyperplanes)} in $V$ is a finite collection of linear hyperplanes in $V$. 
In this article we assume that every arrangement is central unless otherwise specified. 
A \textit{multiplicity} $m$ is a function $m:\A \rightarrow \Z_{>0}$ 
and a pair 
$(\A,m)$ is called a \textit{multiarrangement}. 
Fix a linear form $\alpha_H\ (H \in \A)$ 
in such a way that $\ker(\alpha_H)=H$.
The \textit{logarithmic derivation module} 
$D(\A,m)$ associated with $(\A,m)$ is defined by
$$
D(\A,m):=\{\theta \in \Der_\R(S)\mid\theta(\alpha_H) \in S \cdot 
\alpha_{H}^{m(H)}\ \text{for all}\,\, H \in \A\}.
$$
In general,
 $D(\A,m)$ is not necessarily a free $S$-module.
We say 
that $(\A,m)$ is \textit{free} 
if $D(\A,m)$ is a free $S$-module. 
For a fixed arrangement $\A$, 
a multiplicity $m$ on $\A$ is called \textit{free} 
if a multiarrangement $(\A,m)$ is free. 
Define
\begin{eqnarray*}
{\mathcal{NFM}} (\A)
:=
\{
m : \A \rightarrow \Z_{>0}  \mid m \text{~is~not~a~free~multiplicity}
\}.
\end{eqnarray*} 
The following 
definition was introduced in 
\cite[Definition 5.4]{ATW2}.

\begin{define}
An arrangement $\A$ is called \textit{totally free} 
if 
every multiplicity $m:\A \rightarrow \Z_{>0}$ is a 
free 
multiplicity, or equivalently
${\mathcal{NFM}} (\A)
=
\emptyset$. 
\label{FM}
\end{define}

%

When $\A_{i} $ is an arrangement in $V_{i} \,(i=1, 2)$,
the {\it
product}
 $\A_{1} \times \A_{2} $ is an arrangement
in $V_{1} \oplus V_{2} $ defined 
as in \cite[Definition 2.13]{OT}
by
\[
\A_{1} \times \A_{2} =
\{
H_{1}  \oplus V_{2} \mid H_{1} \in \A_{1} 
\}
\cup
\{
V_{1}  \oplus H_{2} \mid H_{2} \in \A_{2} 
\}.
\]

Our main theorem is as follows:

\begin{theorem}
An arrangement
$\A$ is
totally free if and only if
it has a
decomposition 
$$
\A=\A_1 \times \A_2 \times \cdots \times \A_s,
$$
where each $\A_i$ is an
arrangement in $\K^1$ or $\K^2$.
\label{main}
\end{theorem}

Ziegler showed in \cite[Corollary 7]{Z} that 
$(\A, m)$ is a free multiarrangement whenever
$\ell \le 2.$ 
Note that
 $$D(\A_{1} \times \A_{2}, m) 
\simeq S \cdot D(\A_1,m|_{\A_1}) \oplus 
S \cdot D(\A_2,m|_{\A_2})$$ holds true
as shown in
\cite[Lemma 1.4]{ATW1}. 
Thus
$$
\A_1 \times \A_2 \times \cdots \times \A_s
$$
is known to be totally free if
each $\A_i$ is an arrangement in $\K^1$ or $\K^2$.
Theorem \ref{main} asserts that the converse is
also
 true.
In the next section we will prove
Theorem \ref{main} 
in a stronger form: we will show that $\A$ is decomposed into
one-dimensional arrangements and two-dimensional ones
if 
${\mathcal{NFM}} (\A)
$
is
a finite set. 

Recall that the \textit{intersection lattice} $L(\A)$ is 
the set $\{X=H_1 \cap \cdots \cap H_s
\mid
H_i \in \A,\ s \ge 0\}$ with the reverse inclusion ordering
as in \cite[Definition 2.1]{OT}. 
Then
Theorem \ref{main} implies:
%
%
%
%
\begin{cor}
Whether an arrangement $\A$ 
is totally free or not depends only on
its intersection lattice $L(\A)$. 
\end{cor}

Let $\A$ be a nonempty central arrangement and
$H_{0}\in\A$.  Define
the
{\it deletion} $\A'$ 
and the
{\it restriction} 
$\A''$ 
as in \cite[Definition 1.14]{OT}:
\begin{eqnarray*}
\A' := \A\setminus \{H_{0} \},
\,\,\,\,\,
\A'' := \{
H_{0} \cap H\mid H\in\A'
\}.
\end{eqnarray*}
Because of the characterization
in
Theorem~\ref{main}, 
the total freeness is stable under 
deletion and restriction: 

\begin{cor}
Any subarrangement or restriction of a 
totally free arrangement 
is also totally free.
\end{cor}

A multiarrangement was introduced and studied by Ziegler in \cite{Z}.
The third author proved 
in 
\cite{Y1} and \cite{Y2}
that
the freeness of 
a simple arrangement
is closely
related with the freeness of Ziegler's canonical restriction. 
Recently 
the first and second authors and
 Wakefield
developed
a general theory of free multiarrangements
and introduced
the concept of 
{free
multiplicity}
in \cite{ATW1} and \cite{ATW2}. 
Several
 papers 
including
\cite{A3}, \cite{ANN}, \cite{AY2} and \cite{Y4}
studied the set
of free multiplicities
for a 
fixed arrangement $\A$.
The main theorem (Theorem \ref{main})
in this article
shows
that the set
of free multiplicities
(or ${\mathcal{NFM}} (\A)$)
imposes strong restrictions on
the original
arrangement $\A$.

\medskip

\begin{small} 
\textbf{Acknowledgements}. The first author is supported by the 
JSPS Research Fellowship for Young Scientists. 
The second and third authors have been 
supported in part by Japan Society for the Promotion of Science. 
The authors thank 
Professors
Sergey Yuzvinsky
and
Max Wakefield for helpful
discussions and comments.
We also thank Professor Thomas Zaslavsky for
pointing out an error in an earlier version.
\end{small}

\section{Proof of Theorem~\ref{main} }

First we review a necessary condition for a given multiarrangement to be 
free in Theorem \ref{LMPGMP}. 

Let $(\A,m)$ be a multiarrangement. 
When 
$(\A,m)$ is free, there exists 
a homogeneous basis $\theta_1,\ldots,
\theta_\ell$ for $D(\A,m)$.  The set $\exp(\A,m)$ of 
\textit{exponents} is defined 
by $\exp(\A,m):=(\deg\theta_1,\ldots,\deg \theta_\ell)$, 
where $\deg(\theta_i):=
\deg \theta_i(\alpha)$ for some linear form $\alpha$ 
with $\theta_i(\alpha) \neq 0$. 

Define $L(\A)_2:=\{X \in L(\A)
\mid
\codim_V(X)=2\}$ and $\A_X:=\{H \in \A
\mid
X \subset H\}$. 
For $X \in L(\A)_2$ the multiarrangement $(\A_X,m|_{\A_X})$ is free with exponents 
$(d_1^X,d_2^X,0,\ldots,0)$. Define the 
\textit{second local mixed product} $LMP_2(\A,m)$ 
as in \cite[Definition 4.3]{ATW1} by 
$$
LMP_2(\A,m):=\sum_{X \in L(\A)_2} d_1^X d_2^X.
$$
If $\B$ is a subarrangement of $\A$, then
it is easy to see that
\[
LMP_{2} (\A, m) \ge
LMP_{2} (\B, m|_{\B} ).
\]
Next assume that $(\A,m)$ is free with exponents 
$(d_1,\ldots,d_\ell)$. Define the 
\textit{second global mixed product} $GMP_2(\A,m)$ 
as in \cite[Definition 4.5]{ATW1} by 
$$
GMP_2(\A,m):=\sum_{1 \le i<j \le \ell} d_id_j.
$$

\begin{theorem}
If a multiarrangement $(\A,m)$ is free, then 
$GMP_2(\A,m)=LMP_2(\A,m)$. 
\label{LMPGMP}
\end{theorem}
In fact, Theorem \ref{LMPGMP} is true for any $GMP_k$ and $LMP_k\ (1 \le k \le \ell)$,
see \cite[Corollary 4.6]{ATW1}. 

An arrangement $\A$ is said to be {\em reducible}
if $\A = \A_{1} \times \A_{2} $ for certain
arrangements $\A_{i} $ in $V_{i} \,
(i=1, 2)$.
We say $\A$ is {\it irreducible} if it is not reducible.

\begin{lemma}
\label{tutte} 
Let $\A$ be an irreducible arrangement in $\K^{\ell} $ 
with $\ell\ge 3$.  Then there exists
a subarrangement $\B$  with $|\B| = \ell+1$ 
such that:
\begin{eqnarray*}
\codim_{V} H_{i} \cap 
 H_{j } \cap 
 H_{k } 
 &=&
3
\end{eqnarray*}
whenever
$
H_{i}, H_{j}, H_{k}$ are three hyperplanes in $\B$.
Moreover the arrangement 
$\B$ is not free.
\end{lemma}

\noindent
{\bf Proof.}
Let $H_{0} \in \A$.  Let $\A'$ and $\A''$ be the
deletion and the restriction respectively.
Then either $\A'$ or $\A''$ is irreducible
by Tutte \cite{Tu}
(see also \cite[Theorem 4.3.1]{Ox}).
When $|\A| = \ell+1$, the arrangement
$\A$ itself satisfies the condition.
We will prove by an induction on $|\A|$.
If $\A'$ is irreducible, then $\A'$ contains
$\ell+1$ hyperplanes satisfying the condition.
So we may assume that $\A''$ is irreducible.
Let $\A' = \{H_{1}, H_{2} , \dots , H_{n-1} \}$
and $\overline{H_{i} } := H_{i} \cap H_{0} $ for 
$1\leq i\leq n-1$.  

Suppose $\ell=3$.  
Since $\A''$ is irreducible, $3\leq |\A''|\leq |\A'|$. 

{\it Case 1.}
If $
3\leq
|\A''| = |\A'|$, then 
$
\overline{H_{1} },
\overline{H_{2} },
\dots
,
\overline{H_{n-1} }
 $ are distinct.
Note that
$H_{1} \cap H_{2} \cap\dots\cap H_{n-1} = \{{\mathbf 0}\}$ 
because $\A$ is irreducible.  
Thus we may assume
$H_{1} \cap H_{2} \cap H_{3} =
\{{\mathbf 0}\}.
$ 
So $
\{
H_{0},
H_{1},
H_{2},
H_{3}
\}
 $ 
satisfies the condition.

{\it Case 2.}
If $3\leq |\A''| < |\A'|$, then 
we may assume that
$
\overline{H_{1} },
\overline{H_{2} },
\overline{H_{3} }
 $ 
are distinct and
$
\overline{H_{3} }
=
\overline{H_{4} }$.
Note that $H_{3} \cap H_{4} \subset H_{0}
$. 
Thus
$$
H_{1} \cap H_{2} 
\cap H_{3} 
\cap H_{4} 
=
\overline{H_{1} }
\cap
\overline{H_{2} }
\cap
\overline{H_{3} }
=
\{{\mathbf 0}\}.
$$
Therefore 
either
$
H_{1} \cap H_{2} 
\cap H_{3} 
=
\{{\mathbf 0}\}
$
or
$
H_{1} \cap H_{2} 
\cap H_{4} 
=
\{{\mathbf 0}\}.
$
So
we may
conclude that
either
$
\{
H_{0}, 
H_{1},
H_{2}, 
H_{3} 
\}
$
or
$
\{
H_{0}, 
H_{1},
H_{2}, 
H_{4} 
\}
$
satisfies the condition.

Suppose $\ell \geq 4$. 
Then, by the induction assumption,
there exists a subarrangement
$
\{
\overline{H_{1} },
\overline{H_{2} },
\dots
,
\overline{H_{\ell} }
\}
$
of $
\A''$
satisfying the condition.
Then 
the subarrangement
$
\{
H_{0},
H_{1},
\dots
,
H_{\ell}
\}
$ 
of 
$\A$ satisfies the condition.

Suppose that the arrngement $\B$ is free.  Then the sum of exponents
is equal to $\ell+1$.  Thus one has 
$GMP_2 \leq \binom{\ell}{2} ((\ell+1)/\ell)^{2} $.
We also have $LMP_2 = \binom{\ell+1}{2}  $.
Since $$
LMP_2 - GMP_2 \geq \binom{\ell+1}{2} 
- \binom{\ell}{2} ((\ell+1)/\ell)^{2} = (\ell+1)/2\ell 
> 0,$$ this contradicts Theorem \ref{LMPGMP}.
\owari

\medskip

Recall
$$
{\mathcal{NFM}} (\A)
=
\{
m : \A \rightarrow \Z_{>0}  \mid m \text{~is~not~a~free~multiplicity}
\}.
$$
\begin{prop}
\label{ge3} 
If $\A$ is an 
irreducible arrangement
in $\K^{\ell} $
with $\ell\ge 3$,
then
$
{\mathcal{NFM}} (\A)
$
is an infinite set.
\end{prop}

\noindent
\textbf{Proof}. 
Suppose that $
{\mathcal{NFM}} (\A)
$
is a finite set.
Choose
a subarrangement
$\B
$
of $\A$  with $|\B| = \ell+1$ 
satisfying
the condition in Lemma \ref{tutte}. 
Consider the multiplicity 
$m$ defined by 
\[
m(H)=
\left\{
\begin{array}{rl}
1 & \mbox{if}\ H \not\in \B ,\\
k & \mbox{if}\ H \in \B,
\end{array}
\right.
\]
for every positive integer $k$ .
Since 
$
{\mathcal{NFM}} (\A)
$
is a finite set,
the multiarrangement $(\A,m)$ is free
whenever $k$ is sufficiently large.
Note
$
|
L(\B)_{2} |
=
\binom{\ell+1}{2}.
$
By the definition of 
$LMP_{2} $, 
$$
LMP_2(\A,m) \ge 
LMP_{2}(\B, m|_{\B} ) 
=|L(\B)_{2}| k^{2} = 
\binom{\ell+1}{2} k^{2}.
$$
Let $|\A| = n$.
Then
\[
\sum_{d \in \exp(\A, m)} d =
(k-1)(\ell+1) + n
\]
and
thus
$$
GMP_2(\A,m)
\le
\binom{\ell}{2}
\left\{
\frac{(k-1)(\ell+1)+n}{\ell}  
\right\}^{2} 
=
\frac{(\ell+1)^{2} (\ell-1)}{2 \ell} k^{2}  
+A k + B
$$
with some constants $A$ and $B$.
By Theorem \ref{LMPGMP}
we have
\begin{eqnarray*}
\binom{\ell+1}{2} k^{2}
\le
LMP_2(\A,m)
=
GMP_2(\A,m)
\le
\frac{(\ell+1)^{2} (\ell-1)}{2 \ell} k^{2}  
+A k + B
\end{eqnarray*} 
whenever $k$ is sufficiently large.
This is a contradiction 
because
\[
\binom{\ell+1}{2}
-
\frac{(\ell+1)^{2} (\ell-1)}{2 \ell} 
=
\frac{\ell+1}{2 \ell}  
>0
.
\]
\owari

\medskip

We now prove the following theorem which is
stronger than Theorem \ref{main}.

\begin{theorem}
The following four conditions for a central arrangement $\A$ are equivalent:

(1) $\A$ is totally free, i. e., ${\mathcal{NFM}}(\A)$ is empty,

(2) ${\mathcal{NFM}}(\A)$ is a finite set,

(3) $\A$ has a decomposition 
$$
\A=\A_1 \times \A_2 \times \cdots \times \A_s,
$$
where each $\A_i$ is an
arrangement in $\K^1$ or $\K^2$,

(4) every subarrangement of $\A$ is free.

\end{theorem}

\noindent
\textbf{Proof}. 
The implications 
$(1)\Rightarrow (2)$,
$(3)\Rightarrow (4)$
and
$(3)\Rightarrow (1)$
are obvious.
Thus it is enough to prove that 
$(2) \Rightarrow (3)$
and 
$(4)  \Rightarrow (3).$

$(2) \Rightarrow (3)$:
Suppose that
${\mathcal{NFM}}(\A)$ is a finite set.
Decompose $\A$ into
$$\A_1 \times \A_2 \times \dots\times \A_{s} $$ 
such that each $\A_{i} $ is irreducible.
Since
$$
D(\A,m) \simeq S \cdot D(\A_1,m|_{\A_1}) \oplus 
S \cdot D(\A_2,m|_{\A_2})\oplus\dots\oplus 
S \cdot D(\A_s,m|_{\A_s})$$ holds
by \cite[Lemma 1.4]{ATW1},
each 
$\A_{i} $ is an irreducible arrangement
and
${\mathcal{NFM}}(\A_{i})$ is a finite set.
Thus
Proposition
\ref{ge3} 
shows that 
each arrangement $\A_i$ is in $\K^1$ or $\K^2$. 

$(4) \Rightarrow (3)$:
Decompose $\A$ into irreducible arrangements.
Then each of the irreducible arrangements satisfies 
the assumption (4).
Therefore we may assume that $\A$ is irreducible from
the beginning.
Then, by Lemma \ref{tutte}, 
we may conclude $\ell\le 2$. 
\owari

\end{document}